\documentclass[11pt]{amsart}

\usepackage{amsmath,amsthm,amssymb}
\usepackage[T1]{fontenc}
\usepackage{latexsym}
\usepackage{subfigure, epsfig, epsf}
\usepackage{color,fancyhdr,lastpage,graphicx}
\usepackage{xspace}
\usepackage{wasysym}
\usepackage{setspace}
\usepackage{bbm}

\newcommand{\RR}{\mathbb{R}}

 \newcommand{\mcC}{\mathcal{C}} \newcommand{\mcF}{\mathcal{F}} \newcommand{\mcU}{\mathcal{U}}
  \newcommand{\mcM}{\mathcal{M}} \newcommand{\mcB}{\mathcal{B}}

\newcommand{\ra}{\rightarrow} \newcommand{\la}{\leftarrow}

\newcommand{\st}{such that }

\newcommand{\varep}{\varepsilon}

\renewcommand{\leq}{\leqslant}
\renewcommand{\geq}{\geqslant}
\renewcommand{\la}{\lambda}
\newcommand{\al}{\alpha}
\newcommand{\wrt}{with respect to }
\renewcommand{\st}{such that }

\newcommand{\ssk}{\smallskip }

\newtheorem{thm}{\hspace{-0.15cm} {\sc Theorem} }
\newtheorem{lem}[thm]{\hspace{-0.15cm} {\sc Lemma} }
\newtheorem{prop}[thm]{\hspace{-0.15cm} {\sc Proposition}}

\numberwithin{equation}{section} 

\newenvironment{Dem}{%
    \begin{list}{\hspace{0.6cm}{\sc Proof --}}{%
        \setlength{\topsep}{0pt}%
        \setlength{\leftmargin}{0pt}%
        \setlength{\rightmargin}{0pt}%
        \setlength{\listparindent}{0pt}%
        \setlength{\itemindent}{0pt}%
        \setlength{\parsep}{0pt}%
        \addtolength{\leftmargin}{20pt}%
        \addtolength{\rightmargin}{0pt}%
    } \item }{\hfill{\space $\rhd$}\end{list}\smallskip}

\title{Sensitivity for Smoluchowski equation}
\date{\today}
\author[I.F. Bailleul]{I.F. Bailleul}
\address{Statistical Laboratory, Cambridge, Email: i.bailleul@statslab.cam.ac.uk \\ This research was supported by the EPSRC grant EP/E01772X/1.}
\email{i.bailleul@statslab.cam.ac.uk}
\urladdr{http://www.statslab.cam.ac.uk/~ismael/}


\begin{document}

\begin{abstract}
This article investigates the question of sensitivity of the solutions $\mu_t^\la$ of Smoluchowski equation on $\RR_+^*$ \wrt parameters $\la$ in the interaction kernel $K^{\la}(.,.)$. It is proved that $\mu_t^{\la}$ is a $\mathcal{C}^1$ function of $(t,\la)$ with values in a good space of measures under the hypotheses $K^{\la}(x,y)\leq \varphi(x)\,\varphi(y)$, for some sub-linear function $\varphi$, and $\displaystyle{\int \varphi^{4+\varep}(x)\,\mu_0(dx) < \infty}$, and that the derivative is the unique solution of a related equation.
\end{abstract}

\maketitle

\section{Introduction} 

\noindent \textbf{a) Smoluchowski equation.} Many chemical reactions, such as soot formation \cite{MarkusPatterson} or flame synthesis of organic or inorganic nanoparticles \cite{MarkusNano}, have in common a microscopic mechanism where particles of different masses evolve in a homogeneous medium. Each of them performs a free thermal motion, with diffusivity depending on its mass, until it approaches enough any other particle.  These two particles will then coagulate to create a new one, whose mass will be the sum of the masses of each of its ancestors.

The experimentor have only access to macroscopic quantities such as the concentration of the different masses along time. How can he describe the evolution of these quantities from this microscopic description of the dynamics? Mathematically, we can describe these concentrations as measures $\mu_t$ on the space $\RR^*_+ := (0,+\infty)$  of masses of species. What comes out from experimental measurements are quantities such like the concentration of particles with a mass between such and such number, or, more generally, quantities of the form $(f,\mu_t) \equiv \int f(x)\mu_t(dx)$, for some functions $f$. Smoluchowski has proposed in \cite{Smoluchowski} to describe the evolution of the observations $(f,\mu_t)$ in a well mixed system using some symmetric kernel $K(x,y)$ describing the rates at which coagulations occur:
\begin{equation}
\label{WeakSmo}
\frac{d}{ds}\bigl(f,\mu_s\bigr) = \frac{1}{2}\int\bigl\{f(x+y)-f(x)-f(y)\bigr\}K(x,y)\mu_s(dx)\mu_s(dy).
\end{equation}
Roughly speaking, a particle of mass $x$ coagulates with a particle of mass $y$ at rate $K(x,y)$ to create a particle of mass $x+y$.

\medskip

\noindent \textbf{b) Sensitivity.} The parameters of an experiment are incorporated into the model dynamics \eqref{WeakSmo} as parameters $\la\in\RR^d$ in the interaction kernel $K(\cdot,\cdot) = K^{\la}(\cdot,\cdot)$. Binder granulation a priori requires for instance around $10$ parameters to describe it (\cite{BinderGranulation}). Finding the relevant parameters, given the experimental data (the so-called ``inverse problem'') is the fundamental step which will allow future simulations to provide law cost predictions. Let denote by $\la$ a generic multi-dimensional parameter, $K^\la$ the corresponding coagulation kernel and $\mu_t^\la$ the solution to Smoluchowksi equation associated with $K^\la$. A simple and largely used method for tunning the parameter to data consists in formally applying a method of steepest descent so as to minimize some distance between $\mu_T^\la$ and $\mu_T^{\textrm{obs}}$, in the typical case where we are interested in the value at time $T$ of the system. The measure $\mu_T^{\textrm{obs}}$ is given by experiments. To be effective, the algorithm requires the knowledge of the differential $\sigma^t_\la$ of $\mu_t^\la$ \wrt $\la$ so as to choose the steepest descent direction at each step. Note that $\sigma_t^\la$ is a priori a signed measure. Engineers usually estimate it by a finite difference corresponding to two close values of $\la$. The main approach to do that consists in approximating the differences $\frac{\mu_t^{\la+\epsilon e_i}-\mu_t^\la}{\epsilon}$ (for a basis vector $e_i$ of $\RR^d$) by the corresponding difference for the approximating particle systems -- see \cite{PeterJamesMarkus} for a non-trivial and efficient way of doing that. However, no justification that $\partial_\la\mu_t^\la$ exists has ever been given up to now, which puts the previous investigations on a somewhat hazy mathematical framework. 

\ssk

The aim of this article is to prove that $\mu_t^\la$ is a $\mcC^1$ function of $(t,\la)$ (under proper conditions and in a suitable sense) and that it is the unique solution to some equation (``sensitivity equation''). Not only does this fact put the existing approaches on a firm ground, but it also leads to a new particle approximation \cite{MLSensitivity} which happens to be more accurate than any other method. In the same way as one can associate some finite interacting particle systems to Smoluchowski equation, the so-called Marcus-Lushnikov processes \cite{Marcus}, one can associate a pair of coupled interacting particle systems to the equation associated with the sensitivity, \st their difference converges weakly to a solution of the ``sensitivity equation'', as a consequence of a kind of law of large numbers. The well-posedness of this equation justifies theoretically the use of that particle system for simulating the sensitivity.

\medskip

\noindent \textbf{Notation.} Given a locally bounded non-negative kernel $F(x,y)$ on $\RR^*_+\times\RR^*_+$ and a Radon measures $\mu, \nu$ on $\RR^*_+$, one defines a signed Radon measure $F(\mu,\nu)$ setting
\begin{equation}
\label{DefnKMuMu}
F(\mu,\nu) = \int \bigl\{\delta_{x+y}-\delta_x-\delta_y\bigr\}\,F(x,y)\,\mu(dx)\,\nu(dy).
\end{equation}

\ssk

\noindent

\medskip

\noindent \textbf{c) Strategy for studying the sensitivity of Smoluchowski equation.} We describe in the remainder of this section the approach we use to prove the above mentionned differentiability result. From a mathematical point of view, the main difficulty in solving Smoluchowski equation comes from the fact that whilst the weak formulation \eqref{WeakSmo} is always a well-defined problem (although it may have no solution), it is not easy to find a Banach or a Fr\'echet space of (signed) measures where the differential equation 
\begin{equation}
\label{StrongSmo}
\dot{\mu}_s = \frac{1}{2} K(\mu_s,\mu_s)
\end{equation}
itself is meaningful. This difficulty disappears for bounded kernels, where Smoluchowski equation can be solved in the Banach framework of Radon signed measures equipped with total variation norm. The computation of $\partial_\la\mu_t^{\la}$ is formally straightforward and leads to a representation formula involving essentially only $\{\mu^{\la}_s\}_{s\leq t}$. The map $t\ra\mu_t^{\la}$ solving equation \eqref{StrongSmo}, its derivative with respect to $\la$ solves formally the equation
\begin{equation}
\label{SensitivityEq}
\dot{\sigma}_t^{\la} = K^{\la}\left(\mu^{\la}_t,\sigma_t^{\la}\right) + \frac{1}{2}\partial_{\la}K^{\la}\left(\mu^{\la}_t,\mu^{\la}_t\right) 
\end{equation}
obtained by differentiation of equation \eqref{StrongSmo} with respect to $\la$; we have written $\partial_{\la}K^{\la}(x,y)$ for the partial derivative of $K^{\la}(x,y)$ \wrt $\la$. This equation can be solved, considering first the linearized problem
\begin{equation}
\label{SmoLinLin}
\dot{\rho}_s^{\la} = K^{\la}(\mu^{\la}_s,\rho_s^{\la})
\end{equation}
\noindent before using the variation of constants method.

\smallskip

\textbf{(i)} We introduce a dual evolution equation on functions to study the linear equation \eqref{SmoLinLin}. To that end, define some time dependent operators $\Lambda_s^{\la}$ on functions setting
\begin{equation}
\label{OpEvolutionFunctions}
\Lambda_s^{\la} f(x) = \int \bigl\{f(x+y)-f(x)-f(y)\bigr\}K^{\la}(x,y)\mu_s^{\la}(dy).
\end{equation}
These operators satisfy the identity
$$
\bigl(\Lambda_s^{\la} f,\rho\bigr) = \bigl(f,K(\mu_s^{\la},\rho)\bigr).
$$
Now, if one considers the backward linear equation
$$
\dot{f}_s = -\Lambda_s^{\la} f_s, \quad s\in [0,t] \textrm{ and } f_t=f, 
$$
its solution $\{f_s\}_{0\leq s\leq t}$ depends linearly on $f$, so we can write it in the form $U_{s,t}^{\la}f$, for a linear operator $U_{s,t}^{\la}$. This function $U_{s,t}^{\la}f$ has two important properties. As a function of $t$ it satisfies the identity $\frac{d}{dt} U_{s,t}^{\la}f = U_{s,t}^{\la}\Lambda_t^{\la}f$, and if $\{\rho_s^\la\}_{s\geq 0}$ denotes a solution of equation \eqref{SmoLinLin}, then 
\begin{equation*}
\begin{split}
\frac{d}{ds}\bigl(U_{s,t}^{\la}f,\rho_s^\la\bigr) &= \bigl(-\Lambda_s^{\la}U_{s,t}^{\la}f,\rho_s^{\la}\bigr) + \bigl(U_{s,t}^{\la}f,\dot{\rho}_s^{\la}\bigr) \\
&= -\bigl(U_{s,t}^{\la}f,K(\mu_s^{\la},\rho_s^\la)\bigr) + \bigl(U_{s,t}^{\la}f,K(\mu_s^{\la},\rho_s^\la)\bigr) \\
&= 0.
\end{split}
\end{equation*}
\noindent So we see that the solution to the linear equation \eqref{SmoLinLin} needs to be given by the formula\begin{equation}
\label{PropFundEvolFunction}
\bigl(f,\rho_t^{\la}\bigr) = \bigl(U^{\la}_{0,t}f,\rho_0\bigr).
\end{equation}

\smallskip

\textbf{(ii)} To implement the variation of constants method and solve the affine equation \eqref{SensitivityEq}, introduce as in equation \eqref{OpEvolutionFunctions} the operator
$$
\Lambda_s^{\partial\la} f(x) = \int \bigl\{f(x+y)-f(x)-f(y)\bigr\}\partial_{\la}K^{\la}(x,y)\mu_s^{\la}(dy).
$$
Note the relations
$$
\bigl(\Lambda_s^{\la}f,\mu_s^{\la}\bigr) = \Bigl(f,K^{\la}\bigl(\mu_s^{\la},\mu_s^{\la}\bigr)\Bigr)\quad \textrm{ and } \quad \bigl(\Lambda_s^{\partial\la}f,\mu_s^{\la}\bigr) = \Bigl(f,\partial_{\la}K^{\la}\bigl(\mu_s^{\la},\mu_s^{\la}\bigr)\Bigr).
$$
Defining the measures $\sigma_t^{\la}$ by the formula
\begin{equation}
\label{ReprFormula}
\bigl(f,\sigma_t^{\la}\bigr) = \frac{1}{2}\int_0^t \left(\Lambda_s^{\partial\la}U_{s,t}^{\la}f,\mu_s^{\la}\right)ds
\end{equation}
one sees that it satisfies a weak form of equation \eqref{SensitivityEq}:
\begin{equation*}
\begin{split}
\frac{d}{dt}\bigl(f,\sigma_t^{\la}\bigr) &= \frac{d}{dt}\left(\frac{1}{2}\int_0^t U_{0,s}^{\la}\Lambda_s^{\partial\la}U_{s,t}^{\la}f\,ds\,,\,\mu_0\right) \\
&= \left(\frac{1}{2}\int_0^t U^{0,s}_{\la}\Lambda_s^{\partial\la}U_{s,t}^{\la}\Lambda_t^{\la}f\,ds,\mu_0\right) + \frac{1}{2}\left(U_{0,t}^{\la}\Lambda_t^{\partial\la}f,\mu_0\right) \\
&=  \left(\Lambda_t^{\la}f,\sigma_t^{\la}\right) + \frac{1}{2}\left(\Lambda_t^{\partial\la}f,\mu_t^{\la}\right) \\
&= \Bigl(f,K^{\la}\left(\mu^{\la}_t,\sigma_t^{\la}\right)\Bigr) + \Bigl(f,\frac{1}{2}K^{\partial\la}\left(\mu^{\la}_t,\mu^{\la}_t\right)\Bigr).
\end{split}
\end{equation*}

\medskip

\noindent \textbf{d) Organisation of the article.} How far from full justification is this argument? In the case of uniformly bounded kernels $K^{\la}$, we shall see in section \ref{SensitivityBdedKer} that everything is meaningful in the Banach framework of signed measures equipped with total variation distance, its dual space being the space of bounded functions equipped with the supremum norm. Yet, no such satisfactory framework is available for unbounded kernels; we shall thus use an approximation procedure in section \ref{SensitivityUnboundedKer} to extend the result. The main result (theorem \ref{ThmSensitivity2}) states that the function $(t,\la) \mapsto \mu_t^{\la}$ is a $\mathcal{C}^1$ function with values in a good space of measures and that it is the only solution of a weaker version of equation \eqref{SensitivityEq} under proper conditions.

\smallskip

The idea to investigate a linearized Smoluchowski equation was first used in Kolokoltsov's paper \cite{Kolokoltsov2} to see how $\mu_t$ depends on its initial value. We use here the same tools (theorems \ref{KoloLemma}, \ref{ThmPropagatorsKol1}, \ref{ThmPropagatorsKol2}) as in that paper. We shall compare in section \ref{SectionComments}, \textbf{a)} the present work with the work of Kolokoltsov. Note that the simplified proof of a useful lemma of Kolokoltsov (theorem \ref{KoloLemma}), given in section \ref{SectionComments}, \textbf{b)} and used in section \ref{ConvergenceMeasures}, might be of some interest for itself.

\medskip

\noindent \textbf{Notations.} All functions and measures are defined on $\RR^*_+$ throughout the text.

$\bullet$ We shall use the notation $\mu^{\otimes 2}(dxdy)$ for the product measure $\mu(dx)\mu(dy)$.

$\bullet$ As the expression $f(x+y)-f(x)-f(y)$ will appear numerous times in the text, it will be useful to abbreviate it into $\{f\}(x,y)$. In these terms, the weak version \eqref{WeakSmo} of Smoluchowski equation may be written 
$$
\frac{d}{dt}(f,\mu_t) = \frac{1}{2}\int\{f\}(x,y)K(x,y)\mu_t(dx)\mu_t(dy).
$$

\section{Sensitivity for bounded kernels}
\label{SensitivityBdedKer}

We consider in this section Smoluchowski equation \eqref{WeakSmo} for a family $\{K^\la\}_{\la}$ of interaction kernels, bounded some constant $M$. We recall in section \ref{Section1} why the strong version \eqref{StrongSmo} of Smoluchowski equation is well defined in a good Banach framework. The classical tools of differential equations will then give us for free existence, uniqueness and regularity results of the solutions $\{\mu_t^\la\}_{t\geq 0}$ to equation \eqref{StrongSmo}. We shall then take profit in section \ref{SubSectionSensitivity} of the fact that the derivative $\sigma_t^\la = \partial_\la \mu_t^\la$ solves a time-non-homogeneous \textit{affine} equation to get an explicit formula for it which will be useful in the sequel.

\subsection{Existence and uniqueness in the bounded case: a quick overview}
\label{Section1}

Denote by $B_0$ the Banach space of bounded measurable functions, equipped with the supremum norm $\|.\|_0$. Denote also by $\|\rho\|_0$ the total variation of a signed Radon measure $\rho$, and by 
$$
\mathcal{M}_0 = \{\mu \textrm{ Radon measure }; \|\rho\|_0<\infty\}.
$$ 
Note that $\|\rho\|_0 = \sup\bigl\{(f,\rho)\,;\,f\in B_0, \|f\|_0\leq 1\bigr\}$, and that the space $\left(\mathcal{M}_0,\|.\|_0\right)$ is complete since it is the dual space of the complete space $\bigl(\mcC_b(\RR_+,\RR),\|.\|_{\infty}\bigr)$. We shall denote by $\mathcal{M}_0^+$ the cone of non-negative elements of $\mathcal{M}_0$.

\ssk

\noindent The main reason why everything works well in the bounded case is that as we have $\bigl|\bigl(f,K(\mu,\mu)\bigr)\bigr|\leq 3 \|f\|_0M\|\mu\|_0^2$, for any $f\in B_0$, the Radon measure $K(\mu,\mu)$ belongs to $\mathcal{M}_0$ if $\mu$ does; so Smoluchowski equation \eqref{StrongSmo}: $\dot\mu_s = \frac{1}{2}K(\mu_s,\mu_s)$, is a well defined ordinary differential equation in the Banach space $\mathcal{M}_0$.

\begin{prop}
\label{ExistenceUniquenessBdedKer}
Equation \eqref{StrongSmo} has a well defined flow of solutions in $\bigl(\mathcal{M}_0,\|.\|_0\bigr)$, which preserves the cone $\mathcal{M}_0^+$. The solution $\mu_t$ is defined for all times if $\mu\in\mcM_0^+$.
\end{prop}

\begin{Dem}
{ It suffices to see that the vector field $K$ is locally Lipschitz. But given $\mu$ and $\nu$ in $\mathcal{M}_0$, one can write 
$$
\bigl(\mu^{\otimes 2} - \nu^{\otimes 2}\bigr)(dxdy) = \mu(dx)(\mu-\nu)(dy) + \nu(dy)(\mu-\nu)(dx).
$$
Nothing more is needed to get, for any $f\in B_0$, the inequality
\begin{equation*}
\begin{split}
\bigl|\bigl(f,K(\mu,\mu)\bigr)-\bigl(f,K(\nu,\nu)\bigr)\bigr| &= \Bigl|\int \{f\}(x,y)K(x,y)\bigl(\mu^{\otimes 2} - \nu^{\otimes 2}\bigr)(dxdy)\Bigr| \\
&\leq 3\|f\|_0M\bigl(\|\mu\|_0+\|\nu\|_0\bigr)\|\mu-\nu\|_0,
\end{split}
\end{equation*}
which implies
\begin{equation}
\label{KLip}
\bigl\|K(\mu,\mu)-K(\nu,\nu)\bigr\|_0 \leq 3M\bigl(\|\mu\|_0+\|\nu\|_0\bigr)\|\mu-\nu\|_0.
\end{equation}
To see that $\mu_t$ is non-negative if $\mu_0$ is non-negative we find a non-negative function $\theta_t$ on $\RR^*_+$ such that the transformed measure $\rho_t:=\theta_t\mu_t$ solves a differential equation which preserves $\mathcal{M}_0^+$ in a obvious way\footnote{Which is not the case of Smoluchowski equation. One uses the same method in the study of Boltzmann equation.}. See \cite{James1}, proposition $2.2$, for instance.

Given an initial condition $\mu_0$, denote by $\bigl[0,T(\mu_0)\bigr)$ the maximal interval on which the solution started from $\mu_0$ is defined. If $\mu_0$ is non-negative, one has
$$
\frac{d}{dt}\|\mu_t\|_0 = \frac{d}{dt}(1,\mu_t) = -\frac{1}{2}\int K(x,y)\mu_t(dx)\mu_t(dy) \leq 0
$$
and the path $\{\mu_t\}_{0\leq t<T(\mu_0)}$ stays in a ball where the vector field $K$ is (globally) Lipschitz. This explains why the solution is actually defined on $[0,\infty)$.
}
\end{Dem}

\subsection{Sensitivity}
\label{SubSectionSensitivity}

We prove in this section that if the coagulation kernel depends nicely on a parameter $\la$ then the solution to Smoluchowski equation is a $\mcC^1$ function of $(t,\la)$. It's derivative \wrt $\la$ has a representation involving only $(\mu_s)_{s\geq 0}$.

\subsubsection{Dependence on a parameter}

Let now $\{K^{\la}(.,.)\}_{\la \in\mathcal{U}}$ be a family of symmetric non-negative kernels on $\RR^+_*$ depending in a $\mathcal{C}^2$ way in a parameter $\lambda$ belonging to some open set $\mathcal{U}$ of some $\RR^p$. Denote by $K^{\partial\la}(x,y)$ the derivative of $K^{\la}(x,y)$ with respect to $\la$ and define the Radon signed measure $K^{\partial\la}(\mu,\mu)$ setting 
$$
\bigl(f,K^{\partial\la}(\mu,\mu)\bigr) = \int\{f(x+y)-f(x)-f(y)\}K^{\partial\la}(x,y)\,\mu(dx)\mu(dy).
$$
Denote by $\bigl[0,T^{\la}(\mu_0)\bigr)$ the maximal interval on which the solution to Smoluchowski equation \eqref{StrongSmo} with interaction kernel $K^{\la}(\cdot,\cdot)$ started from $\mu_0$ is defined

\begin{thm}[Sensitivity for bounded kernels]
\label{ThmSensitivity}
Suppose $K^{\la}(\cdot,\cdot)$ and its first two derviatives are bounded by a constant $M$, uniformly in $\la\in\mcU$. Then the map $(t,\la)\in\bigl[0,T^{\la}(\mu_0)\bigr)\times\mcU\mapsto\mu_t^{\la}\in\left(\mathcal{M}_0,\|.\|_0\right)$ is differentiable with respect to $\la$ and its derivatives $\sigma_t^{\la}$ (called ``sensitivity'') is the unique solution of the equation 
\begin{equation}
\label{SensitivityEq2}
\dot{\sigma}_t^{\la} = K^{\la}\left(\mu^{\la}_t,\sigma_t^{\la}\right) + \frac{1}{2}K^{\partial\la}\left(\mu^{\la}_t,\mu^{\la}_t\right).
\end{equation}
\end{thm}

\begin{Dem}
As is classically done in the study of ordinary differential equations in Banach spaces (e.g. consult \cite{RHMartin}), the result is a consequence the following four properties.

\begin{enumerate}
   \item For each $\mu\in\mathcal{M}_0$, the map $\la\in\mathcal{U}\mapsto K^{\la}(\mu,\mu)\in \left(\mathcal{M}_0,\|.\|_0\right)$ is differentiable, with a derivative $K^{\partial\la}(\mu,\mu)\in (\mathcal{M}_0,\|.\|_0)$ depending continuously on $\mu\in (\mathcal{M}_0,\|.\|_0)$.
   \item The map $(s,\la)\mapsto \mu_s^{\la}\in\left(\mathcal{M}_0,\|.\|_0\right)$ is continuous on $[0,T]\times\mathcal{U}$.
   \item The linear map $\nu\mapsto K(\mu_s,\nu)$ takes $\left(\mathcal{M}_0,\|.\|_0\right)$ into itself and has a uniformly bounded norm for $s\in [0,T]$. The same result holds for the map $\nu\mapsto K^{\partial\la}(\mu_s,\nu)$.
   \item Let $C$ be a compact set of $\left(\mathcal{M}_0,\|.\|_0\right)$. There exists an $\left(\mathcal{M}_0,\|.\|_0\right)$-valued function $O_2(\mu,\mu')$ such that $\|O_2(\mu,\mu')\|_0\leq m\|\mu-\mu'\|_0^2$ for some constant $m$, and 
\begin{equation}
\label{RegK}
\forall\,\mu,\,\mu'\in C, \quad K^{\la_0}(\mu',\mu')-K^{\la_0}(\mu,\mu) = 2K^{\la_0}(\mu,\mu'-\mu) + O_2(\mu,\mu').
\end{equation}
$K^{\partial\la_0}$ has the same property.
\end{enumerate}

We prove points 1 and 2 and leave the elementary proofs of points 3 and 4 to the reader. 

\noindent 1. Given $f\in B_0$, apply Taylor formula in a small neighbourhood $\mathcal{V}$ of $\la_0$ to get
\begin{equation*}
\begin{split}
&\left|\bigl(f,K^{\la}(\mu,\mu)-K^{\la_0}(\mu,\mu)-(\la-\la_0)K^{\partial\la}(\mu,\mu)\bigr)\right| =\\ 
& \left|\int \{f\}(x,y)\bigl(K^{\la}(x,y)-K^{\la_0}(x,y)-(\la-\la_0)K^{\partial\la}(x,y)\bigr)\mu(dx)\mu(dy)\right| \\
&\leq 3\|f\|_0\frac{|\la-\la_0|^2}{2}\,\underset{\widetilde{\la}\in\mathcal{V}}{\max}\bigl|\partial_{\widetilde{\la}}^2K^{\widetilde{\la}}(x,y)\bigr|\|\mu\|_0^2.
\end{split}
\end{equation*}
This proves the differentiability assertion. The map $\mu\in\mcM_0\ra K^{\partial\la}(\mu,\mu)$ can be seen to be locally Lipschitz using the same reasonning as was used in the proof of proposition \ref{ExistenceUniquenessBdedKer} to prove that the vector field $K$ is locally Lipschitz.

\noindent 2. It is a classical fact in dynamics\footnote{Consult Martin's book \cite{RHMartin} for instance.} that it is sufficient to check that the map $(\la,\mu)\in\mathcal{U}\times\mathcal{M}_0 \rightarrow K^{\la}(\mu,\mu)$ is locally Lipschitz to get the continuity of $(s,\la)\mapsto \mu_s^{\la}\in\left(\mathcal{M}_0,\|.\|_0\right)$. Writing 
$$
K^{\la}(\mu,\mu) - K^{\la'}(\nu,\nu) = K^{\la}(\mu,\mu) - K^{\la}(\nu,\nu) + \left(K^{\la}-K^{\la'}\right)(\nu,\nu),
$$ 
and using inequality \eqref{KLip}, Taylor formula, and the fact that $\underset{x,y\,;\,\ell}{\sup} \bigl|K^{\partial\ell}(x,y)\bigr| \leq M$, one obtains
$$
\|K^{\la}(\mu,\mu) - K^{\la'}(\nu,\nu)\|_0 \leq 3M\bigl(\|\mu\|_0+\|\nu\|_0\bigr)\|\mu-\nu\|_0 + 3M\|\nu\|_0^2 |\la-\la'|.
$$
\end{Dem}

\subsubsection{A representation formula for the sensitivity}

We fix $\mu_0$ throughout this section and work on a fixed time interval $[0,T]\subset\bigl[0,T^{\la}(\mu_0)\bigr)$, for all $\la\in\mcU$. As explained in the introduction, one can solve explicitly equation \eqref{SensitivityEq2} solving first its linearized verion before using the variation of constant method. The first step is made solving a dual problem to the homogeneous equation, on the space $B_0$.

\ssk

\noindent \textbf{a) Dual linearized Smoluchowski equation.} Define for each $\la\in\mathcal{U}$, a time-dependent linear vector field $\Lambda^{\la}_s$ on $B_0$, setting for any $f\in B_0$
\begin{equation}
\label{Op}
\Lambda^{\la}_sf(x) = \int \bigl\{f(x+y)-f(x)-f(y)\bigr\}K^{\la}(x,y)\mu^{\la}_s(dy).
\end{equation}
As $\|\Lambda_s^{\la}\|_0 \leq 3M(1,\mu^{\la}_s)\leq 3M\|\mu_0\|_0$, and $\mu_s^\la$ depends continuously on $s$, the vector field $\Lambda^{\la}_s$ on $B_0$ is continuous \wrt $f\in B_0$ and $s$. So, given some time $t>0$, the backward and forward differential equations 
\begin{equation}
\label{EqGenerator}
\dot{f}_s(x) = -\Lambda^{\la}_s f_s\,(x), \quad f_t \textrm{ given,}  
\end{equation}
are meaningful in $B_0$, and elementary results on linear differential equations on Banach spaces give the following proposition\footnote{Consult Martin's book \cite{RHMartin}.}.

\begin{prop}
\label{Propagator}
The differential backwards and forwards equations \eqref{EqGenerator} in $\left(B_0,\|.\|_0\right)$ have a unique solution, defined for all time. It is of the form $f_s = U_{s,t}^{\la}f_t$, for a continuous linear operator $U_{s,t}^{\la}$ on $B_0$, with norm $\leq e^{3M\|\mu_0\|_0 |t-s|}$. We also have for any $f\in B_0$
\begin{equation}
\label{DepenUWrtT}
\frac{d}{dt}U_{s,t}^{\la}f = U_{s,t}^{\la}\Lambda_t^{\la}f.
\end{equation}
\end{prop}

\noindent This operator $U^{\la}_{s,t}$ can be used to solve explicitly the linear equation on $\mcM_0$
$$
\dot{\rho}_s^{\la} = K^{\la}(\mu^{\la}_s,\rho_s^{\la});
$$
this equation has a unique solution on the time interval $[0,T]$ as the time non-homogeneous vector field  $K^{\la}(\mu^{\la}_s,\cdot)$ is continuous and bounded. Indeed, one gets from Smoluchowski equation \eqref{StrongSmo} and equation \eqref{EqGenerator}
\begin{equation*}
\begin{split}
& \frac{d}{ds}\bigl(U^{\la}_{s,t}f,\rho_s\bigr) = -\bigl(\Lambda_s^{\la}U^{\la}_{s,t}f,\rho_s\bigr) + \bigl(U^{\la}_{s,t}f,\dot{\rho}_s\bigr) \\
& = -\bigl(U^{\la}_{s,t}f,K(\mu_s,\rho_s)\bigr) + \bigl(U^{\la}_{s,t}f,K(\mu_s,\rho_s)\bigr) \\
& = 0;
\end{split}
\end{equation*}
so the identity $\bigl(U_{0,t}^{\la}f,\rho_0\bigr) = (f,\rho_t^{\la})$ holds for any $f\in B_0$; thus
$$
\rho_t^{\la} = \bigl(U_{0,t}^{\la}\bigr)^*\rho_0.
$$

\noindent \textbf{b) A representation formula for $\sigma_t^{\la}$.} The second step to solve the affine equation \eqref{SensitivityEq2} is to use the variation of constant method as explained in the introduction. The following lemma will be used in the way.

\begin{lem}
\label{PetitLem}
The function $t\in [0,T] \mapsto \sigma_t^{\la}\in \mathcal{M}_0$ is the only solution in $\bigl(\mathcal{M}_0,\|.\|_0\bigr)$ of the weak differential equation
$$
\forall\,f\in B_0, \quad \frac{d}{dt}(f,\sigma_t) = \Bigl(f,K^{\la}\bigl(\mu_t^{\la},\sigma_t\bigr)\Bigr) + \frac{1}{2}\Bigl(f,K^{\la}(\mu_t^{\la},\mu_t^{\la}\bigr)\Bigr), \quad \sigma_0 \textrm{ given}.
$$
\end{lem}

\begin{Dem}
{ Note first that since the function $t\in [0,T] \mapsto \sigma_t^{\la}\in \mathcal{M}_0$ satisfies the strong equation \eqref{SensitivityEq2} it also satisfies the above weak equation. Given two solutions $\sigma_t$ and $\overline\sigma_t$ of the latter, one has for any $f\in B_0$
$$
\bigl(f,\sigma_t-\overline\sigma_t\bigr) = \int_0^t \Bigl(f,K^{\la}\bigl(\mu_s^{\la},\overline\sigma_s-\sigma_s\bigr)\Bigr)ds = \int_0^t \bigl(\Lambda_s^{\la}f,\overline\sigma_s-\sigma_s\bigr)ds.
$$
But as the operator $\Lambda_s^{\la}$ on $\bigl(B_0,\|.\|_0\bigr)$ has norm $\leq 3M\|\mu_0\|_0$, we must have
$$
\bigl(f,\sigma_t-\overline\sigma_t\bigr) \leq 3M\|\mu_0\|_0\|f\|_0 \int_0^t \|\sigma_t-\overline\sigma_t\|_0ds,
$$
and so
$$
\|\sigma_t-\overline\sigma_t\|_0 \leq 3M\|\mu_0\|_0 \int_0^t \|\sigma_s-\overline\sigma_s\|_0\,ds.
$$
One deduces from Gronwall's formula that $\overline\sigma_t=\sigma_t$.
}
\end{Dem}

\noindent Define the map $\Lambda^{\partial\la}_s$ on $B_0$ by the formula
$$
\Lambda^{\partial\la}_sf(x) = \int \bigl\{f(x+y)-f(x)-f(y)\bigr\}K^{\partial\la}(x,y)\mu^{\la}_s(dy);
$$
notice that the identities
$$
\bigl(\Lambda^{\partial\la}_sf,\mu^{\la}_s\bigr) = \bigl(f,K^{\partial\la}_s(\mu^{\la}_s,\mu^{\la}_s)\bigr),\quad\textrm{and}\quad \bigl(\Lambda^{\la}_sf,\mu^{\la}_s\bigr) = \bigl(f,K^{\la}_s(\mu^{\la}_s,\mu^{\la}_s)\bigr), \quad f\in B_0.
$$

\begin{prop}[Representation formula for the sensitivity]
\label{PropExplicitSensitivity}
One has 
\begin{equation}
\label{ExplicitSensitivity}
(f,\sigma_t^{\la}) = \frac{1}{2}\int_0^t \left(\Lambda_s^{\partial\la}U_{s,t}^{\la}f,\mu_s^{\la}\right)ds
\end{equation}
for any $f\in B_0$.
\end{prop}

\begin{Dem}
{ Denote temporarily by $\widehat{\sigma}_t^{\la}$ the measure $f\in B_0 \mapsto \frac{1}{2}\int_0^t \left(\Lambda_s^{\partial\la}U_{s,t}^{\la}f,\mu_s^{\la}\right)ds$; it belongs to $\mathcal{M}_0$. The following calculus is fully justified in the Banach framework of $\left(B_0,\|.\|_0\right)$. For any $f\in B_0$, one has
\begin{equation*}
\begin{split}
 \frac{d}{dt} \bigl(f,\widehat{\sigma}_t^{\la}\bigr) &= \frac{d}{dt}\left(\frac{1}{2}\int_0^t \Lambda_s^{\partial\la}U_{s,t}^{\la}f\,ds\,,\,\mu_s^{\la}\right)= \left(\frac{1}{2}\int_0^t \Lambda_s^{\partial\la}U_{s,t}^{\la}\Lambda_t^{\la}f\,ds,\mu_s^{\la}\right) + \frac{1}{2}\left(\Lambda_t^{\partial\la}f,\mu_t^{\la}\right) \\
&=  \left(\Lambda_t^{\la}f,\widehat{\sigma}_t^{\la}\right) + \frac{1}{2}\left(\Lambda_t^{\partial\la}f,\mu_t^{\la}\right) \\
&= \Bigl(f,K^{\la}\left(\mu^{\la}_t,\widehat{\sigma}_t^{\la}\right)\Bigr) + \Bigl(f,\frac{1}{2}K^{\partial\la}\left(\mu^{\la}_t,\mu^{\la}_t\right)\Bigr).\end{split}
\end{equation*}
Since $\widehat{\sigma}_t^{\la}$ satisfies a weak version of equation \eqref{SensitivityEq2} it coincides with $\sigma_t^{\la}$ according to lemma \ref{PetitLem}.
}
\end{Dem}

\section{From bounded to unbounded kernels}
\label{SensitivityUnboundedKer}

We shall now drop the boundedness hypothesis on the kernels $K^\la$. Yet, to get some control on the interaction rates, we shall make the hypothesis that one has
 \begin{equation}
\label{Condition1}
\forall\,\la\in\mcU, \forall\,x,y\in\RR^*_+, \quad K^{\la}(x,y)\leq\varphi(x)\,\varphi(y)
\end{equation}
for some sub-additive function $\varphi$(\footnote{We have $\varphi(x+y)\leq \varphi(x)+\varphi(y)$, for all $x,y\in\RR^*_+$.}), greater than $1$. We shall also suppose that 
\begin{equation}
\label{Condition2}
K^{\partial\la}(x,y)\leq\varphi(x)\,\varphi(y).
\end{equation}
Last, we shall suppose the existence of a (small) $\varep>0$ \st
\begin{equation}
\label{ConditionMoment}
(\varphi^{4+\varep},\mu_0) < \infty.
\end{equation}

In his paper \cite{James1}, J. Norris proved that $(\varphi^2,\mu^{\la}_t)$ remains finite on some time interval $\bigl[0,T(\mu_0)\bigr)$ if $(\varphi^2,\mu_0)$ is finite. The same argument shows that $(\varphi^{4+\varep},\mu^{\la}_t)$ also remains finite (on a possibly different time interval, still denoted $\bigl[0,T(\mu_0)\bigr)$) if $(\varphi^{4+\varep},\mu_0)$ is finite. Given some $T<T(\mu_0)$ denote by $C(T)$ a positive constant such that
\begin{equation}
\label{MajorMoments}
\forall\,t\leq T, \quad  (\varphi^{4+\varep},\mu_t) \leq C(T).
\end{equation}
The function $\varphi$ being greater than $1$, the other moments $(\varphi^p,\mu^{\la}_t)$, with $1\leq p\leq 4+\varep$, are also bounded above by $C(T)$ on $[0,T]$.

\medskip

\noindent In order to estimate the tail behaviour of measures, we introduce the following spaces of measures, indexed by non-negative reals $p$: 
$$
\mathcal{M}_p = \bigl\{\mu\,;\,\|\mu\|_p:=\bigl(\varphi^p,|\mu|\bigr)<\infty\bigr\}.
$$
Using this notation condition \eqref{MajorMoments} reads: $\mu_t\in\mathcal{M}_{4+\varep}\subset \mathcal{M}_1$, for all $0\leq t \leq T$. To compare the behaviour of non-bounded functions with the behaviour of $\varphi$, one defines the increasing family of function spaces, indexed by non-negative reals $p$:
$$
B_p = \Bigl\{f\,;\, \sup\frac{|f|}{\varphi^p}<\infty\Bigr\};
$$
we shall write $\|f\|_p$ for this supremum. Note that $\|\mu\|_p = \sup\{(f,\mu)\,;\,f\in B_p,\,\|f\|_p\leq 1\}$. The purpose of this section is to prove our main result.

\begin{thm}[Sensitivity for unbounded kernels]
\label{ThmSensitivity2}
Assume conditions \eqref{Condition1}, \eqref{Condition2} and the moment condition \eqref{ConditionMoment}. Then the map $(t,\la)\in [0,T]\times\mathcal{U}\mapsto \mu_t^{\la}\in\bigl(\mathcal{M}_1,\,\|.\|_1\bigr)$, is a $\mathcal{C}^1$ function and its derivative $\sigma_t^\la$ satisfies the following equation for any $f\in B_0$. {\small
$$
\left(f,\sigma_t^{\la}\right) = \left(f,\sigma_0^{\la}\right) + \int_0^t \int \{f\}(x,y)K^{\la}(x,y)\mu^{\la}_s(dx)\sigma_s^{\la}(dy)ds + \frac{1}{2}\int_0^t\int \{f\}(x,y)K^{\partial\la}(x,y)\mu^{\la}_s(dx)\mu_s^{\la}(dy)ds
$$}
The function $\sigma_\cdot^\la$ is the only $\bigl(\mathcal{M}_1,\,\|.\|_1\bigr)$-valued solution of this equation.
\end{thm}

\noindent This statement will be proved by an approximation procedure. Let $\left\{K^{\la\,;\,N}\right\}_{N\geq 0}$ be a sequence of bounded symmetric kernels converging towards $K$, and such that $\partial_{\la}K^{\la\,;\,N}$ and $\partial^2_{\la}K^{\la\,;\,N}$ are also bounded, with $\left|K^{\la\,;\,N}(x,y)\right|$ and $\left|\partial_{\la}K^{\la\,;\,N}(x,y)\right|$ bounded above by $\varphi(x)\,\varphi(y)$. Let $\mu_t^{\la\,;\,N}$ and $\sigma_t^{\la\,;\,N}$ be the measures associated with $K^{\la\,;\,N}$ and $\partial_{\la}K^{\la\,;\,N}$, constructed in section \ref{SensitivityBdedKer}. Theorem \ref{ThmSensitivity2} will be proved by showing that
\begin{enumerate}
   \item the map $(t,\la)\in [0,T]\times\mathcal{U} \mapsto \mu_t^{\la\,;\,N}\in\bigl(\mathcal{M}_1,\|.\|_1\bigr)$ is, for each $N$, a $\mathcal{C}^1$ function, and $\partial_{\la}\,\mu_t^{\la\,;\,N} = \sigma_t^{\la\,;\,N}$ in $\bigl(\mathcal{M}_1,\|.\|_1\bigr)$.
   \item the sequence $\bigl\{\mu_t^{\la\,;\,N}\bigr\}_{N\geq 0}$ converges towards $\mu_t^{\la}$ in $\bigl(\mathcal{M}_1,\|.\|_1\bigr)$, uniformly with respect to $(t,\la)\in [0,T]\times\mathcal{U}$;
   \item the sequence $\bigl\{\sigma_t^{\la\,;\,N}\bigr\}_{N\geq 0}$ of its derivatives converges in $\bigl(\mathcal{M}_1,\|.\|_1\bigr)$ towards some $\sigma_t^{\la}$, uniformly with respect to $(t,\la)\in [0,T]\times\mathcal{U}$. 
\end{enumerate}
Points 2 and 3 will be proved sections \ref{ConvergenceMeasures} and \ref{ConvergencePropagators} respectively. We prove the first point here. Denote by $M$ an upper bound of $K^{\la\,;\,N}$. Notice first that the inequality $\bigl|\{f\}(x,y)\bigr|\leq 2\|f\|_1\bigl(\varphi(x)+\varphi(y)\bigr)$, gives for any $\mu\in\mcM_1$
\begin{equation}
\begin{split}
\bigl|\bigl(f,K^{\la\,;\,N}(\mu,\mu)\bigr)\bigr| &\leq 2M\|f\|_1\int\bigl(\varphi(x)+\varphi(y)\bigr)\mu(dx)\mu(dy) \\
&\leq 4M \|f\|_1\|\mu\|_1^2;
\end{split}
\end{equation}
so the Radon measure $K^{\la\,;\,N}(\mu,\mu)$ belongs to $\mathcal{M}_1$ if $\mu$ does. Now, the following inequalities enable us to see that the vector field $\mu\mapsto K^{\la\,;\,N}(\mu,\mu)$ on $\bigl(\mathcal{M}_1,\|.\|_1\bigr)$ is Lipschitz. The function $f\in B_1$ has norm no greater than $1$ and $\mu,\nu \in\mathcal{M}_1$.
\begin{equation*}
\begin{split}
\Bigl|\Bigl(f,K^{\la\,;\,N}(\mu,\mu)-&K^{\la\,;\,N}(\nu,\nu)\bigr)\Bigr)\Bigr| = 2M\int \bigl(\varphi(x)+\varphi(y)\bigr)\bigl(|\mu|(dx)|\mu-\nu|(dy) + |\nu|(dy)|\mu-\nu|(dx)\bigr) \\
&\leq 2M\bigl(\|\mu\|_1\|\mu-\nu\|_0 + \|\mu\|_0\|\mu-\nu\|_1 + \|\nu\|_0\|\mu-\nu\|_1 + \|\nu\|_1\|\mu-\nu\|_0\bigr) \\
&\leq 4M \bigl(\|\mu\|_1 + \|\nu\|_1\bigr)\|\mu-\nu\|_1.
\end{split}
\end{equation*}
The differentiability of the map $\la\in\mathcal{U} \mapsto \mu_t^{\la\,;\,N}\in\bigl(\mathcal{M}_1,\|.\|_1\bigr)$ can be proved in the same way as was done in section \ref{SubSectionSensitivity} in the framework of $\bigl(\mathcal{M}_0,\|.\|_0\bigr)$. To prove the continuity of $\mu_t^{\la\,;\,N}$ and $\sigma_t^{\la\,;\,N}$ with respect to $(t,\la)\in [0,T]\times\mathcal{U}$, one checks that the vector fields appearing in equations \eqref{StrongSmo} and \eqref{SensitivityEq2} are Lipschitz in $(\la,\mu)\in\mcU\times\mcM_1$, mimicking what was done in the proof of theorm \ref{ThmSensitivity} in the framework of $\mcU\times\mcM_0$. This completes the proof of the first point.

\ssk

Note that the operators $\Lambda_s^{\la\,;\,N}$ and $\Lambda_s^{\partial\la\,;\,N}$ are bounded in $\bigl(\mathcal{M}_1,\|.\|_1\bigr)$, with norm no greater than $4M\|\mu_s^{\la}\|_1$, so that the representation formula for $\sigma_t^{\la\,;\,N}$ given in \eqref{ExplicitSensitivity} also holds in $\bigl(\mathcal{M}_1,\|.\|_1\bigr)$. The remainder of this section is dedicated to the proofs of points 2 and 3. After a preliminary result in section \ref{SubSectionProp}, we prove a stronger version of point 2, useful in the sequel. The proof of point 3 is made in section \ref{ConvergencePropagators}.

\smallskip

As we shall prove these results for a fixed $\la$, we shall drop the $\la$ in $\mu_t^{\la}$ and $\sigma_t^{\la}$ in the sequel. The following elementary result will be used repeatedly; its proof is left to the reader.

\begin{lem}
For any $p\geq 1$ and any $f\in B_p$, $\bigl|\{f\}(x,y)\bigr|\leq 2^p\|f\|_p\bigl(\varphi^p(x)+\varphi^p(y)\bigr).$
\end{lem}
\noindent As a last remark, note that the measures $\mu_t^N$ satisfy for any $0\leq t\leq T$ and $N\geq 0$ the same moment inequality \eqref{MajorMoments} as $\mu_t$.

\subsection{Convergence of $\mu_t^N$ to $\mu_t$ in $\left(\mathcal{M}_{2+\varep}, \|.\|_{2+\varep}\right)$}
\label{SubSectionProp}\label{ConvergenceMeasures}

Let $\{\mu_t\}_{0\leq t< T(\mu_0)}$ be the solution given by Norris' theorem; choose $T<T(\mu_0)$. It is worth noting that using dominated convergence and the moment estimate \eqref{MajorMoments}, the measures $\{\mu_t\}_{0\leq t\leq T}$ satisfy the weak version \eqref{WeakSmo} of Smoluchowski equation for any $f\in B_{3+\varep}$. We start this section showing that they depend regularly on $t$.



\begin{prop}
The path $\{\mu_t\}_{0\leq t \leq T}$ is a $\mathcal{C}^1$ path in $\bigl(\mathcal{M}_{2+\varep},\|.\|_{2+\varep}\bigr)$.
\end{prop}

\begin{Dem}
{ One proves that the path $\{\mu_t\}_{0\leq t \leq T}$ is \textit{\textbf{1)}} Lipschitz in $\bigl(\mathcal{M}_{3+\varep},\|.\|_{3+\varep}\bigr)$, \textit{\textbf{2)}} $\mathcal{C}^1$ in $\bigl(\mathcal{M}_{2+\varep},\|.\|_{2+\varep}\bigr)$.

\smallskip

\noindent \textit{\textbf{1)}} Take a function $f\in B_{3+\varep}$. One establishes the following inequalities using the inequality $K(x,y) \leq \varphi(x)\varphi(y)$ and the sub-additivity of $\varphi$. 
\begin{equation*}
\begin{split}
\bigl|(f,\mu_t-\mu_s)\bigr| &\leq \frac{1}{2}\int_s^t \int \bigl|\{f\}(x,y)\bigr|K(x,y)\mu_r(dx)\mu_r(dy)\,dr \\
&\leq \frac{c_{\varep}\|f\|_{3+\varep}}{2} \int_s^t \int \bigl\{\varphi^{3+\varep}(x)+\varphi^{3+\varep}(y)\bigr\}\varphi(x)\varphi(y)\mu_r(dx)\mu_r(dy)\,dr \\
&\leq 2c_{\varep}\|f\|_{3+\varep} \int_s^t \int \varphi^{4+\varep}(x)\varphi(y)\mu_r(dx)\mu_r(dy)\,dr \\
&\leq 2c_{\varep}\|f\|_{3+\varep} (\varphi,\mu_0)\underset{s\leq r\leq t}{\sup} \|\mu_r\|_{4+\varep}\,|t-s|.
\end{split}
\end{equation*}
Taking the supremum of the left hand side, with $\|f\|_{3+\varep}\leq 1$, this shows that the path $\{\mu_t\}_{0\leq t \leq T}$ is Lipschitz in $\bigl(\mathcal{M}_{3+\varep},\|.\|_{3+\varep}\bigr)$, with Lipschitz constant $\leq 2c_{\varep}C(T)^2$.

\smallskip

It follows from this fact that the formula
$$
(f,\nu_t) := \frac{1}{2}\int\{f\}(x,y)K(x,y)\mu_t(dx)\mu_t(dy)
$$
defines an element $\nu_t$ of $\bigl(\mathcal{M}_{2+\varep},\|.\|_{2+\varep}\bigr)$ which is continuous with respect to $t$. Indeed, since one has for any $f\in B_{2+\varep}$, 
\begin{equation*}
\begin{split}
\bigl|(f,\nu_t&-\nu_s)\bigr| = \frac{1}{2}\left|\int \{f\}(x,y)K(x,y) \bigl\{\mu_t(dx)(\mu_t-\mu_s)(dy) + \mu_s(dy)(\mu_t-\mu_s)(dx)\bigr\}\right| \\
&\leq \frac{c'_{\varep}\|f\|_{2+\varep}}{2} \int \bigl(\varphi^{2+\varep}(x)+\varphi^{2+\varep}(y)\bigr)\varphi(x)\varphi(y) \bigl(\mu_t(dx)|\mu_t-\mu_s|(dy) + \mu_s(dy)|\mu_t-\mu_s|(dx)\bigr) \\ 
&\leq 2c'_{\varep}\|f\|_{2+\varep}C(T) \|\mu_t-\mu_s\|_{3+\varep},
\end{split}
\end{equation*}
we have $\|\nu_t-\nu_s\|_{2+\varep} \leq 8c_{\varep}c'_{\varep}C(T)^3\,|t-s|$. 

\smallskip

\textit{\textbf{2)}} Finally, write for any $f\in B_{2+\varep}$
$$
\bigl(f,\mu_t-\mu_s-(t-s)\nu_s\bigr) = \int_s^t (f,\nu_r-\nu_s)dr,
$$
and note that the integral is uniformly $o(t-s)$, for $\|f\|_{2+\varep}\leq 1$; this proves that the path $\{\mu_t\}_{0\leq t\leq T}$ is differentiable, as a path in $\bigl(\mathcal{M}_{2+\varep},\|.\|_{2+\varep}\bigr)$, with continuous derivative $\nu_t$.
}
\end{Dem}

\noindent We shall use this result in the form: \textit{The path $\{\varphi^{2+\varep}\mu_t\}_{0\leq t \leq T}$ is a $\mathcal{C}^1$ path in $\bigl(\mathcal{M}_0,\|.\|_0\bigr)$}. This enables us to apply a useful lemma of Kolokoltsov (see the appendix of \cite{Kolokoltsov1}) of which we give a clear and short proof in section \ref{SectionComments}. 

\begin{lem}[Kolokoltsov \cite{Kolokoltsov1}]
\label{KoloLemma}
Let $\{\rho_s\}_{0\leq s \leq T}$ be a $\mathcal{C}^1$ path in $\left(\mathcal{M}_0,\|.\|_0\right)$, with derivative $\{\dot{\rho}_s\}_{0\leq s\leq T}$. There exists a $\{\pm 1,0\}$-valued measurable function $(s,x)\in\RR_+\times\RR^*_+\mapsto\varep_s(x)$ such that we have 
\begin{itemize}
   \item $\|\rho_t\|_0 = \|\rho_0\|_0 + \int_0^t(\varep_s,\dot{\rho}_s)\,ds,\quad \textrm{for any } t\in[0,T]$,
   \item $\bigl(f,|\rho_t|\bigr) = (f\varep_t,\rho_t)$, for all $f\in\mcB, t\in [0,T]$.
\end{itemize}
\end{lem}

\begin{prop}
\label{CvgceMuN}
The sequence of measures $\{\mu_t^N\}_{N\geq 0}$ converges to $\mu_t$ in $\bigl(\mathcal{M}_{2+\varep},\|.\|_{2+\varep}\bigr)$, uniformly with respect to $t\in [0,T]$.  
\end{prop}

\begin{Dem}
Applying Kolokoltsov's lemma to the $\mathcal{C}^1$ path $\bigl\{\varphi^{2+\varep}(\mu_t^N-\mu_t)\bigr\}_{0\leq t \leq T}$ in $\bigl(\mathcal{M}_0,\|.\|_0\bigr)$, and denoting by $\varep_s^N$ the function given by theorem \ref{KoloLemma}, we can write
\begin{equation*}
\begin{split} 
\|\mu_t^N-\mu_t\|_{2+\varep}  &= \int\varphi^{2+\varep}(x)|\mu_t^N-\mu_t|(dx) = \int_0^t \bigl(\varep_s^N\varphi^{2+\varep},\dot{\mu}^N_s-\dot{\mu}_s\bigr)ds \\
&= \int_0^t\{\varep_s^N\varphi^{2+\varep}\}(x,y)\left(K^N(x,y){\mu^N_s}^{\otimes 2} - K(x,y)\mu_s^{\otimes 2}\right)(dx\otimes dy) \\
&= \int_0^t\int \{\varep_s^N\varphi^{2+\varep}\}(x,y)K^N(x,y)\bigl({\mu^N_s}^{\otimes 2}-\mu_s^{\otimes 2}\bigr)(dx\otimes dy)\,ds \\
   &+ \int_0^t\int \{\varep_s^N\varphi^{2+\varep}\}(x,y)\bigl(K^N-K\bigr)(x,y)\mu_s^{\otimes 2}(dx\otimes dy)\,ds.
\end{split}
\end{equation*}
The second term converges to $0$ by dominated convergence and the fact that $\|\mu_s\|_{3+\varep}$ is bounded; call it $o_N(1)$. To handle the first term, write it as
\begin{equation*}
\begin{split}
&  \int_0^t\int \{\varep_s^N\varphi^{2+\varep}\}(x,y)K^N(x,y)\Bigl((\mu_s^N-\mu_s)(dx)\mu_s^N(dy) + \mu_s(dx)(\mu_s^N-\mu_s)(dy)\Bigr) \\
&= \int_0^t\int \{\varep_s^N\varphi^{2+\varep}\}(x,y)K^N(x,y)\varep_s^N(x)|\mu_s^N-\mu_s|(dx)\bigl(\mu_s+\mu_s^N\bigr)(dy) =: (*); \\
\end{split}
\end{equation*}
we have used the symmetry of the expressions with respect to $x$ and $y$. Now, using the fact that $\bigl|\varep_s^N\bigr|\leq 1$, one can find some constant $C_{\varep}$ such that 
\begin{equation*}
\begin{split}
\{\varep_s^N\varphi^{2+\varep}\}(x,y)\varep^N_s(x) &\leq \varep^N_s(x)\varphi^{2+\varep}(x+y)-\varphi^{2+\varep}(x)-\varep_s^N(y)\varep_s^N(x)\varphi^{2+\varep}(y) \\
                               &\leq \varphi^{2+\varep}(x+y)-\varphi^{2+\varep}(x)-\varep_s^N(y)\varep_s^N(x)\varphi^{2+\varep}(y). 
\end{split}
\end{equation*}
To deal with the upper bound, note that there exists a constant $C_{\varep}$ \st the inequality 
$$
(a+b)^{\al}-a^{\al} \leq C_\al\left(a^{\al-1}b+b^{\al}\right).
$$ 
holds for any $a,b \geq 0$. It follows that 
\begin{equation*}
\begin{split}
\{\varep_s^N\varphi^{2+\varep}\}(x,y)\varep^N_s(x) \leq C_{\varep}\bigl(\varphi^{2+\varep}(y)+\varphi^{1+\varep}(x)\varphi(y)\bigr),
\end{split}
\end{equation*}
so
\begin{equation*}
\begin{split}
(*) &\leq c_{\varep}\int_0^t\int \bigl(\varphi^{2+\varep}(y)+\varphi^{1+\varep}(x)\varphi(y)\bigr)K^N(x,y)\bigl(\mu_s+\mu_s^N\bigr)(dy)|\mu_s^N-\mu_s|(dx)\,ds \\
&\leq c_{\varep}\int_0^t \Bigl(2\bigl(\|\mu_s\|_{3+\varep}\vee\|\mu_s^N\|_{3+\varep}\bigr)\|\mu_s^N-\mu_s\|_1 + 2\bigl(\|\mu_s\|_2\vee\|\mu_s^N\|_2\bigr)\|\mu_s^N-\mu_s\|_{2+\varep}\Bigr)\,ds \\
&\leq 4C_{\varep}C(T)\int_0^t \|\mu_s^N-\mu_s\|_{2+\varep}\,ds.
\end{split}
\end{equation*}
Putting the pieces together, we have obtained
$$
\|\mu_t^N-\mu_t\|_{2+\varep}  \leq o_N(1) + 4C_{\varep}C(T)\int_0^t \|\mu_s^N-\mu_s\|_{2+\varep}\,ds,
$$ 
where $o_N(1)$ is uniform in $t\in [0,T]$; Gronwall's lemma enables to conclude.
\end{Dem}

All the estimates above do not depend on the implicit parameter $\la$; this proposition proves (a stronger version of) point 2 in our strategy of proof.

\subsection{Convergence of $\sigma_t^N$ to $\sigma_t$ in $\bigl(\mathcal{M}_1, \|.\|_1\bigr)$}
\label{ConvergencePropagators}

We prove the third point of our strategy in this section. For that purpose, we rely crucially on the representation formula \eqref{ExplicitSensitivity} for $\sigma_t$ for bounded kernels, as it brings back the problem of proving the convergence of $\sigma_t^N$ to a convergence problem for $(\mu_s^N)_{0\leq s\leq t}$ and its functionals $U_{s,t}^N$. 
Given $\ell\geq 0$, denote by $B^0_{\ell}$ the set of real-valued functions $f$ on $\RR_+$ \st $\frac{|f|}{\varphi^{\ell}}$ is bounded and converges to $0$ at infinity.

\begin{prop}
\label{ApplThmKol}
\begin{enumerate}
   \item There exists a uniformly bounded family of operators $\{U_{s,t}\}_{0\leq s\leq t\leq T}$ on $\bigl(B_3^0,\|.\|_3\bigr)$ \st the functions $s, t \mapsto U_{s,t}f$ are differentiable in $\bigl(B_3^0,\|.\|_0\bigr)$, when $f\in B_{1+\varep}$, with derivatives $-\Lambda_sU_{s,t}f$ and $U_{s,t}\Lambda_tf$, respectively.
   \item These operators $U_{s,t}$ preserve $B_{1+\varep}^0$, and are bounded in $\bigl(B_{1+\varep}^0,\|.\|_{1+\varep}\bigr)$.
\end{enumerate}
\end{prop}

\begin{Dem}
This proposition is a direct application of theorems \ref{ThmPropagatorsKol1} and \ref{ThmPropagatorsKol2} on propagators, in section \ref{Appendix}; we apply them to the two pairs $\bigl(\varphi^{1+\varep},\varphi^3\bigr)$ and $\bigl(\varphi^{\frac{1}{2}},\varphi^{1+\varep}\bigr)$. We adopt the notations 
$$
\emph{{\bf J}}_s f(x) \equiv \int\bigl\{f(x+y)-f(x)\bigr\}K(x,y)\mu_s(dy), \quad \emph{{\bf M}}_s f(x) \equiv \int f(y)K(x,y)\mu_s(dy)
$$
used in section \ref{Appendix}.

\ssk

\textit{1.} Applying theorems \ref{ThmPropagatorsKol1} and \ref{ThmPropagatorsKol2}, we only need to check that the inequalities
\begin{itemize} 
   \item $\emph{{\bf J}}_s\varphi^{\al} \leq C(\al)\|\mu_s\|_{\al+1}\varphi^{\al}$,
   \item $\bigl|\emph{{\bf M}}_s\bigl(\varphi^{\al}\bigr)\bigr|\leq \|\mu_s\|_{\al+1}\varphi$,
   \item for any $f\in B_{\beta}$, $\quad\quad \emph{{\bf J}}_sf \leq 2^{\beta+1}\bigl(\varphi^{\beta+1}(x)\|\mu_s\|_1+\varphi(x)\|\mu_s\|_{\beta+1}\bigr)$.
\end{itemize}
 hold for any $\al$ and $\beta \geq 1$, which is done by elementary algebra.

\smallskip

\noindent \textit{2.} To apply theorems \ref{ThmPropagatorsKol1} and \ref{ThmPropagatorsKol2} to the pair  $\bigl(\varphi^{\frac{1}{2}},\varphi^{1+\varep}\bigr)$, one needs to verify that $\displaystyle{\emph{{\bf J}}_s\varphi^{\frac{1}{2}} \leq \frac{C(T)}{2}\varphi^{\frac{1}{2}}}$. This can be done by writing
\begin{equation*}
\begin{split}
\int \bigl\{\varphi^{\frac{1}{2}}(x+y)-\varphi^{\frac{1}{2}}(x)\bigr\}K(x,y)\mu_s(dy) &\leq \int\left\{{\bigl(\varphi(x)+\varphi(y)\bigr)}^{\frac{1}{2}}-\varphi^{\frac{1}{2}}(x)\right\}K(x,y)\mu_s(dy) \\
&\leq \int \frac{\varphi(y)}{2\varphi^{\frac{1}{2}}(x)}\varphi(x)\varphi(y)\mu_s(dy) = \frac{\|\mu_s\|_{2}}{2}\varphi^{\frac{1}{2}}(x) \\
&\leq \frac{C(T)}{2}\varphi^{\frac{1}{2}}(x) .
\end{split}
\end{equation*}
\end{Dem}

\noindent Theorem \ref{ThmPropagatorsKol2} provides us with an additional information: $U_{s,t}$ and all its approximations $U_{s,t}^N$ have a norm on $B^0_{1+\varep}$ controlled by the right hand side of equation \eqref{UpperBound}, which is independent of $N$.

\ssk

Since $U_{s,t}$ sends $B_{1+\varep}^0$ in itself, and $\Lambda_s^{\partial\la}$ is easily verified to be a bounded operator from $B_{1+\varep}$ into $B_{2+\varep}$, with a uniformly bounded norm for $0\leq s\leq t\leq T$, the formula  
\begin{equation}
\label{FormulaSigma}
(f,\sigma_t) = \frac{1}{2}\int_0^t \left(\Lambda_s^{\partial\la}U_{s,t}f,\mu_s\right)ds
\end{equation}
defines a measure $\sigma_t$ belonging to $\mathcal{M}_1$. By proposition \ref{ApplThmKol}, the quantities $\bigl(f,\sigma_t^N\bigr)$ and $\bigl(f,\sigma_t\bigr)$ are bounded uniformly in $t\in [0,T]$, $N\geq 0$ and $\la\in\mathcal{U}$, given any $f\in B_1$.

\begin{thm}
\label{ThmRepere}
The sequence $\bigl\{\sigma_t^N\bigr\}_{N\geq 0}$ converges to $\sigma_t$ in $\left(\mathcal{M}_1, \|.\|_1\right)$,  uniformly for $t\in [0,T]$. 
\end{thm}

\begin{Dem}
We need to prove that the limit 
\begin{equation*}
\bigl(f,\sigma_t^N\bigr) =  \frac{1}{2}\int_0^t\int \left\{U_{s,t}^Nf\right\}K^{\partial\la\,;\,N}\bigl(\mu_s^N,\mu_s^N\bigr)\,ds \underset{N,+\infty}{\longrightarrow}  \frac{1}{2}\int_0^t\int \left\{U_{s,t}f\right\}K^{\partial\la}\bigl(\mu_s,\mu_s\bigr)\,ds = (f,\sigma_t) 
\end{equation*}
holds uniformly for $\|f\|_1\leq 1$ and $0\leq t\leq T$. If one can prove that $U_{s,t}^Nf$ converges to $U_{s,t}f$ in $B_{1+\varep}$, uniformly in $0\leq s\leq t\leq T$, then 
\begin{itemize}
   \item the inequality $\big|K^{\partial\la\,;\,N}\big|(x,y)\leq \varphi(x)\varphi(y)$, 
   \item and the fact that $\mu_s^N$ converges to $\mu_s$ in $\bigl(\mathcal{M}_{2+\varep}, \|.\|_{2+\varep}\bigr)$, uniformly in $s\in [0,T]$, 
\end{itemize}
will enable us to apply dominated convergence to get the result. We are thus led to prove that there exists a decreasing sequence $\{a_N\}_{N\geq 0}$, converging to $0$, such that one has 
$$
\bigl\|U_{s,t}f-U_{s,t}^Nf\|_{1+\varep} \leq a_N\bigr\|f\|_1,
$$
for any $0\leq s\leq t\leq T$ and any $f\in B_1$.

Since $f\in B_1 \subset B_{1+\varep}$ one can use the differentiability property of $U_{s,t}$ as a function of $s$ and $t$ to write
$$
U_{s,t}f-U_{s,t}^Nf = \int_s^t \frac{d}{du}\left(U_{s,u}U_{u,t}^N\right)f\,du = \int_s^t \left(U_{s,u}\bigl(\Lambda_u-\Lambda_u^N\bigr)U_{u,t}^N\right)f\,du.
$$
As $U_{u,t}^Nf$ belongs to $B^0_{1+\varep}$, with a norm uniformly controlled for $\|f\|_1\leq 1$, and as $U_{s,u}$ is a uniformly bounded operator on $B_{2+\varep}$, it suffices to prove that \textit{there exists a decreasing sequence $\{a_N\}_{N\geq 0}$ converging to $0$ such that one has 
$$
\bigl\|\bigl(\Lambda_u-\Lambda_u^N\bigr)g\bigr\|_{2+\varep} \leq a_N,
$$
for any $g\in B_{1+\varep}$, with $\|g\|_{1+\varep}\leq 1$.} To prove this fact, write
\begin{equation*}
\begin{split}
&\Bigl|\bigl(\Lambda^{\la}_u-\Lambda^{\la\,;\,N}_u\bigr)g(x)\Bigr| = \left|\int\{g\}(x,y)\Bigl(K(x,y)\mu_s(dy)-K^N(x,y)\mu^N_s(dy)\Bigr)\right| \\
&\leq c_{\varep}\int\bigl(\varphi^{1+\varep}(x)+\varphi(y)\bigr)\Bigl(K(x,y)|\mu_s-\mu_s^N|(dy)+|K-K^N|(x,y)\mu_s(dy)+2\varphi(x)\varphi(y)|\mu_s^N-\mu_s|(dy)\Bigr) \\
&\leq c_{\varep}\varphi^{2+\varep}(x)\,\|\mu_s^N-\mu_s\|_1 + c_{\varep}\varphi(x)\,\|\mu_s-\mu_s^N\|_{2+\varep} + c_{\varep}\varphi^{1+\varep}(x)\left(|K-K^N|(x,.),\mu_s\right) \\ 
&\quad\quad\quad\quad\quad\quad\quad + c_{\varep}\Bigl(\varphi^{1+\varep}(\cdot)\,\big|K-K^N\big|(x,\cdot),\mu_s\Bigr) + c_{\varep}\varphi^{2+\varep}(x)\,\|\mu_s^N-\mu_s\|_1 + c_{\varep}\varphi(x)\,\|\mu_s^N-\mu_s\|_{2+\varep}.
\end{split}
\end{equation*}
This formula makes it clear that we shall get the existence of these $a_N$'s if we can prove that the sequence of functions $x\mapsto\left(\varphi^{1+\varep}(\cdot)\,\big|K-K^N\big|(x,\cdot),\mu_s\right)$ converges to $0$ in $B_{2+\varep}$ as $N\ra +\infty$. This fact is clearly seen on the following inequality where $M$ is an arbitrary positive constant.
\begin{equation*}
\begin{split}
\frac{1}{\varphi^{2+\varep}(x)}\int\varphi^{1+\varep}(y)|K-K^N|&(x,y)\mu_s(dy) \leq \frac{1}{\varphi^{2+\varep}(x)}\int\varphi^{2+\varep}(y)\varphi(x)\textbf{1}_{\varphi(x)\varphi(y)\geq N}\,\mu_s(dy) \\
&\leq \frac{1}{\varphi^{1+\varep}(x)}\int\varphi^{2+\varep}(y)\textbf{1}_{\varphi(x)\varphi(y)\geq N}\,\mu_s(dy) \\
&\overset{\varphi\geq 1}{\leq} \frac{\|\mu_s\|_{2+\varep}}{M^{1+\varep}}\textbf{1}_{\varphi(x)\geq M} + \left(\int\varphi^{2+\varep}(y)\textbf{1}_{\varphi(y)\geq \frac{N}{M}}\mu_s(dy)\right)\textbf{1}_{\varphi(x)\leq M}
\end{split}
\end{equation*}
\end{Dem}

\noindent Proposition \ref{CvgceMuN} and theorem \ref{ThmRepere} together prove point (2) and (3) of our strategy of proof for theorem \ref{ThmSensitivity2}, showing that $\mu_t^{\la}$ is a $\mcC^1$ function of its arguments. To complete the proof of theorem \ref{ThmSensitivity2}, it remains to prove that $\sigma_t^{\la}$ is the unique solution in $\mcM_1$ of the equation {\small
\begin{equation}
\label{EqDerivee}
\left(f,\sigma_t^{\la}\right) = \left(f,\sigma_0^{\la}\right) + \int_0^t \int \{f\}(x,y)K^{\la}(x,y)\,\mu^{\la}_s(dx)\sigma_s^{\la}(dy)ds + \frac{1}{2}\int_0^t\int \{f\}(x,y)K^{\partial\la}(x,y)\,\mu^{\la}_s(dx)\mu_s^{\la}(dy)ds
\end{equation}}
where $f$ is any \textit{bounded} function. 

We have seen in section \ref{SubSectionSensitivity} that this identity holds if one replaces $\sigma_t^{\la}$ and $\mu_t^{\la}$ by  $\sigma_t^{\la\,;\,N}$ and $\mu_t^{\la\,;\,N}$ resspectively. Use then the above convergence results $\sigma_t^{\la\,;\,N}\ra \sigma_t^{\la}$, in $\mcM_1$, and $\displaystyle{\mu_t^{\la\,;\,N}\ra\mu_t^{\la}}$, in $\mcM_{2+\varep}$, together with the inequalities
\begin{equation*}
\begin{split}
& \Bigl|\bigl(f,K^{\la}(\mu_t^{\la},\sigma_t^{\la})\bigr) - \bigl(f,K^{\la}(\mu_t^{\la\,;\,N},\sigma_t^{\la\,;\,N})\bigr)\Bigr| \leq 3\|f\|_{\infty}\,\Bigl(\|\mu_t^{\la}-\mu_t^{\la\,;\,N}\|_1\,\|\sigma_t^{\la}\|_1 + \|\mu_t^{\la\,;\,N}\|_1\,\|\sigma_t^{\la}-\sigma_t^{\la\,;\,N}\|_1 \Bigr), \\
& \Bigl|\bigl(f,K^{\partial\la}(\mu_t^{\la},\mu_t^{\la})\bigr) - \bigl(f,K^{\partial\la}(\mu_t^{\la\,;\,N},\mu_t^{\la\,;\,N})\bigr)\Bigr| \leq 3C(T)\,\|f\|_{\infty}\,\Bigl(\|\mu_t^{\la}-\mu_t^{\la\,;\,N}\|_1+ \|\sigma_t^{\la}-\sigma_t^{\la\,;\,N}\|_1 \Bigr),
\end{split}
\end{equation*}
to pass to the limit properly. 

\medskip

To prove uniqueness of the solution to equation \eqref{EqDerivee} in $\bigl(\mcM_1,\|\cdot\|_1\bigr)$ it suffices to show that the equation 
$$
\forall\,f\in B_c, \quad (f,\gamma_t) = \int_0^t\int\bigl\{f(x+y)-f(x)-f(y)\bigr\}K(x,y)\mu_s(dx)\gamma_s(dy)ds
$$
has at most one solution in $\bigl(\mcM_1,\|\cdot\|_1\bigr)$. We have written here $B_c$ for the set of bounded Borel functions with compact support. Rewrite this equation under the form
$$
(f,\gamma_t) = \int_0^t\bigl(\Lambda_sf,\gamma_s\bigr)ds.
$$
Repeating the proof of corollary \ref{ApplThmKol}, it is seen that there exists bounded propagators $U_{s,t}$ on $\bigl(B^0_1,\|\cdot\|_1\bigr)$ \st the function $s\in [0,t]\mapsto U_{s,t}f$ solves the equation $\frac{d}{ds}U_{s,t}f = -\Lambda_sU_{s,t}f$ for any $f\in B_c(\subset B^0_1)$ and $t\in [0,T]$. It follows that the expression $(U_{s,t}f,\gamma_s)$ is well defined and that 
$$
\frac{d}{ds}(U_{s,t}f,\gamma_s) = \bigl(-\Lambda_sU_{s,t}f,\gamma_s\bigr) + \bigl(\Lambda_sU_{s,t}f,\gamma_s\bigr) = 0;
$$
so $(f,\gamma_t)=(U_{0,t}f,\gamma_0)$, implying the uniqueness of $\gamma_t$. This ends the proof of theorem \ref{ThmSensitivity2}.

\section{Comments}
\label{SectionComments}

\subsection{Related works.} 

One can see the main roots of theorem \ref{ThmSensitivity2} in section $4$ of Kolokoltsov's pioneering article \cite{Kolokoltsov2} on the central limit theorem for the Marcus-Lushnikov dynamics. He develops in this section tools for the analysis of the rate of convergence of the semi-group of Marcus-Lushnikov process to the semi-group of solutions of Smoluchowski equation. Recall the Marcus-Lushnikov process $\bigl\{X^n_t\bigr\}_{t\geq 0}$ is a strong Markov jump process on the space of discrete measures whose jumps are as follows. If its state at time $t$ is $\displaystyle{\frac{1}{n}\sum\delta_{x_i(t)}}$, for $i$ in a finite set $I_t$ depending on $t$, define, for $i<j$ in $I_t$, independent exponential random times $T_{ij}$ with parameter $\displaystyle{\frac{K\bigl(x_i(t),x_j(t)\bigr)}{n}}$ and set
$$
T = \min \{T_{ij}\,;\,i<j\}.
$$
The process remains constant on the time interval $[t,t+T[$ and has a jump $\displaystyle{\frac{1}{n}\bigl(\delta_{x_p(t)+x_q(t)} - \delta_{x_p(t)}-\delta_{x_q(t)}\bigr)}$ at time $t+T$, if $T=T_{pq}$. The dynamics then starts afresh. The convergence of this sequence $\{X^n\}_{n\geq 0}$ of processes to the deterministic solution of Smoluchowski equation was first proved under general conditions in \cite{James1}. Yet, no fine analysis of the convergence of the corresponding semi-group was done before \cite{Kolokoltsov2}. We explain roughly his idea to see how similar equations to the 'variation' equations \eqref{StrongSmo}, \eqref{SensitivityEq} appear in his context.

Suppose we are in a situation where existence and uniqueness of solutions to Smoluchowski equation hold into a proper sense, and denote by $\{T_t\}$ and $\{T_t^n\}$ the semi-groups of Smoluchowski and Marcus-Lushnikov dynamics. Also, denote by $L$ and $L^n$ their generators. Then, given any (good) function $F$ and a measure $\mu$
$$
\bigl(T_t-T_t^n\bigr)F\,(\mu) = \int_0^t \Bigl(T_{t-s}^n\bigl(L^n-L\bigr)T_sF\Bigr)(\mu)\,ds.
$$
The choice of a function $F$ of the form $F(\mu) = \int g({\bf x})\mu^{\otimes k}(d{\bf x})$, for some symmetric function $g$ of $k$ variables, provides a 'measure' of the moments of $\mu$. One has $T_sF\,(\mu) = F(\mu_s)$, where $\mu_0 = \mu$.

Introducing some derivation operation $\delta$ on functions on measures: 
$$
\delta F(\mu\,;\,x) = \underset{\varep\rightarrow 0}{\lim}\,\frac{F(\mu+\varep\delta_x)-F(\mu)}{\varep},
$$
one can write for any function $G$
\begin{equation}
\label{DiffGenerators}
\bigl(L^n-L\bigr)G\,(\mu) = -\frac{1}{2n}\int \bigl(\delta G(\mu\,;\,2x)- 2\delta G(\mu\,;\,x)\bigr)K(x,x)\mu(dx) + O\bigl(n^{-3/2}\bigr).
\end{equation}
One thus sees that taking $G=T_sF$, with the above $F$, leads to consider the quantity
$$
\delta\bigl((g,\mu_t^{\otimes k})\bigr) = k\,\bigl((g,\mu_t^{\otimes k-1}\otimes\delta\mu_t)\bigr),
$$
where 
$$
\delta\mu_t = \underset{\varep\rightarrow 0}{\lim}\,\frac{\mu_t(\mu+\varep\delta_x)-\mu_t(\mu)}{\varep}
$$
is 'the' derivative of $\mu_t$ \wrt its initial condition. Terms of the form $\delta(\delta\mu_t)$ arise in the $O\bigl(n^{-3/2}\bigr)$ term of equation \eqref{DiffGenerators}. This analysis brings back the estimate of $\bigl(T_t-T_t^n\bigr)F\,(\mu)$ to estimates on $\mu_s, \delta\mu_s$ and $\delta^2\mu_s$. To do so, Kolokoltsov shows that $\delta\mu_s$ is a solution of the linear equation 
$$
\frac{d}{ds}\delta\mu_s = K\bigl(\mu_s,\delta\mu_s\bigr)
$$
in some sense, and that $\delta^2\mu_s$ is a solution of the affine equation
$$
\frac{d}{ds}\delta^2\mu_s = K\bigl(\mu_s,\delta^2\mu_s\bigr) + K\bigl(\delta\mu_s,\delta\mu_s\bigr)
$$
in some sense. The tools used to solve these equations are essentially the same as those used above; the reader may will find the details given here helpful to unzip the section $4$ of \cite{Kolokoltsov2}. We have used yet a slightly different approach in the implementation of the variation of constant method. Note also that we have been able to go from the framework of 'sub-linear' kernels of \cite{Kolokoltsov2}: $K(x,y)\leq C(1+x+y)$, to the framework of an essentially 'sub-multiplicative' kernel: $K(x,y)\leq \varphi(x)\,\varphi(y)$, an improvement which is of some practical interest.

\subsection{Kolokoltsov's lemma.} 

This paragraph contains a simple proof of Kolokoltsov's lemma, which was used in a crucial way to prove a uniqueness result in the original article \cite{Kolokoltsov1} where it was first introduced. We prove it here in a slightly less general framework than in \cite{Kolokoltsov1}, sufficient for our purposes as well as for its use in \cite{Kolokoltsov1}; the gain in clarity and volume of the proof is substantial.

\smallskip

Let $(\Omega,\mcF)$ be a measurable space with a $\sigma$-algebra $\mcF$ generated by a filtration $\{\mcF_n\}_{n \geq 0}$ made up of finite $\sigma$-algebras. We shall denote by $\{A_n^p\}_p$ the atoms of $\mathcal{F}_n$. We shall write $(\mcM,\|.\|)$ for the space of finite signed-measures on $(\Omega,\mcF)$, equipped with the total variation distance. We shall define, for each $n\geq 1$, the total variation of a measure \wrt $\mcF_n$:
$$
\forall\,\mu\in\mcM, \quad \|\mu\|_{(n)} = \sup\bigl\{(f,\mu)\,;\,f\in\mcF_n,\,|f|\leq 1\bigr\}.
$$
These quantities have the property 
\begin{equation}
\label{IneqNorms}
\forall\,\mu\in\mcM, \quad \|\mu\|_{(n)} \underset{n+\infty}{\longrightarrow} \|\mu\|.
\end{equation}

Recall that the topological dual space of $(\mcM,\|.\|)$ is the space $(\mcB,|.|)$ of bounded measurable functions on $(\Omega,\mcF)$, equipped with the supremum norm. We shall write $\widehat{\mcB}$ for the set of bounded functions $g$ on $[0,T]\times \Omega$, with norm $\widehat{\|g\|} = \sup\bigl\{g_s(x)\,;\,s\in [0,T],\,x\in\Omega\bigr\}$, and shall define $\left(\widehat{\mcM},\|.\|_{TV}\right)$ as the space of finite signed measures on $[0,T]\times \Omega$, equipped with the total variation norm.

\begin{thm}[Kolokoltsov's lemma \cite{Kolokoltsov1}, Appendix]
\label{KoloLemma}
Let $\{\rho_s\}_{0\leq s \leq T}$ be a $\mathcal{C}^1$ path in $\left(\mathcal{M},\|.\|\right)$, with derivative $\{\dot{\rho}_s\}_{0\leq s\leq T}$. There exists a $\{\pm 1,0\}$-valued measurable function $\varep_s(x)$ such that we have 
\begin{itemize}
   \item $\|\rho_t\| = \|\rho_0\| + \int_0^t(\varep_s,\dot{\rho}_s)\,ds,\quad \textrm{for any } t\in[0,T]$,
   \item $\forall\,f\in\mcB, \forall\,t\in [0,T],\quad \bigl(f,|\rho_t|\bigr) = (f\varep_t,\rho_t)$.
\end{itemize}
\end{thm}

We shall make use of the following elementary lemma in the course of the proof of theorem \ref{KoloLemma}.

\begin{lem}
\label{LemFunction}
By convention, $\emph{sgn}(0)=0$. We have for any $\mathcal{C}^1$ function $g : \RR_+ \rightarrow \RR$
$$
\bigl|g(t)\bigr| = \bigl|g(0)\bigr| + \int_0^t \emph{sgn}\bigl(g(s)\bigr)\,g'(s)ds.
$$
\end{lem}

\begin{Dem}
{ Using lemma \ref{LemFunction} in each set $A_n^p$, we can define a $\{\pm 1,0\}$-valued function $s\mapsto\varep^{n;\,p}_s$ such that
$$
\bigl|\rho_t(A_n^p)\bigr| = \bigl|\rho_0(A_n^p)\bigr| + \int_0^t \varep^{n;\,p}_s\dot{\rho}_s(A_n^p)\,ds.
$$
Define then the function $\varep_s^n(x)$ as being equal to $\varep^{n;\,p}_s$ on $A_n^p$; the preceding identity yields
\begin{equation}
\label{PreEqKoProof}
\|\rho_t\|_{(n)} = \|\rho_0\|_{(n)} + \int_0^t (\varep^n_s,\dot{\rho}_s)\,ds.
\end{equation}
The functions $\varep^n$ belong to the set $\widehat{\mcB}$ of bounded functions on $[0,T]\times \Omega$, and have supremum norm no greater than $1$. Using the duality between $\widehat{\mcM}$ and $\widehat{\mcB}$ provided by integration, equation \eqref{PreEqKoProof} can be written
\begin{equation}
\label{EqKoProof}
\|\rho_t\|_{(n)} = \|\rho_0\|_{(n)} + \bigl(\varep^n,\dot{\rho}_s\otimes ds\bigr).
\end{equation}
Now, since $\left(\widehat{\mcB},\widehat{\|.\|}\right)$ is the topological dual space of $\left(\widehat{\mcM},\|.\|_{TV}\right)$, its unit sphere is weakly-$^*$ compact. We can thus find a sub-sequence $\{\varep^{n_k}\}_{k\geq 1}$ and an element $\varep$ of $\widehat{\mcB}$, with norm less than $1$, \st
$$
\forall\,\mu\in\widehat{\mcM}, \quad (\varep^{n_k},\mu)\underset{k+\infty}{\longrightarrow} (\varep,\mu).
$$
Together with formulas \eqref{IneqNorms} and \eqref{EqKoProof}, this convergence result, applied to the measures $\dot{\rho}_s(dx)\otimes {\bf 1}_{[0,t]}(s)ds$, gives
\begin{equation}
\label{LemKoloIdentity}
\forall\,t\in [0,T],\quad \|\rho_t\| = \|\rho_0\| + \int_0^t (\varep_s,\dot{\rho}_s)\,ds.
\end{equation}

\smallskip

To prove the second point of theorem \ref{KoloLemma}, remark that since 
$$
\int_0^T \|\rho_s\|_{(n)}ds = \int_0^T (\varep^n_s,\rho_s)ds = (\varep^n,\rho_s\otimes {\bf 1}_{[0,T]}ds),
$$
we have 
$$
 \||\rho_s|\otimes ds\|_{TV} = (\varep,\rho_s\otimes {\bf 1}_{[0,T]}ds).
$$
It follows that 
$$
\varep_s\rho_s = |\rho_s|
$$
for almost all $s$. Define $\varep_s$ to be equal to $\frac{d\rho_s}{d|\rho_s|}$ on the exceptional set. This modification of $\varep_s$ preserves identity \eqref{LemKoloIdentity} and proves the second point of theorem \ref{KoloLemma}.
}
\end{Dem}

\section{Appendix on propagators}
\label{Appendix}

We collect in this appendix the material on propagators needed in section \ref{ConvergencePropagators} to prove the convergence of $\sigma_t^N$ to $\sigma_t$ in $\bigl(\mcM_1,\|\cdot\|_1\bigr)$. Recall that a \textit{\textbf{propagator}} is a family $\{U_{s,t}\}_{s\leq t}$ of operators \st $U_{tt} = \textrm{Id}$ and one has $U_{st}U_{tr} = U_{sr}$ for all $s\leq t\leq r$. Define an (a priori) unbounded operator on functions setting
$$
\Lambda_sf(x) = \int\{f\}(x,y)K(x,y)\mu_s(dy).
$$
Theorem \ref{ThmPropagatorsKol2} below states conditions under which the backward/forward differential equation
\begin{equation}
\label{EvolEq}
\dot{u}_s = -\Lambda_s u_s, \;0\leq s\leq t\leq T, \quad u_t \textrm{ given,}
\end{equation}
can be solved in some Banach space of functions. Some notations are needed. Set 
\begin{equation}
\begin{split}
\emph{{\bf J}}_s f(x) \equiv \int\bigl\{f(x+y)-f(x)\bigr\}K(x,y)\mu_s(dy) &= \int f(x+y)K(x,y)\mu_s(dy) - \left(\int K(x,y)\mu_s(dy)\right)f(x) \\
&\equiv \emph{{\bf L}}_sf(x) - \tau_s(x)f(x), 
\end{split}
\end{equation}
and
\begin{equation}
\emph{{\bf M}}_s f(x) \equiv \int f(y)K(x,y)\mu_s(dy), \quad T_s(x) \equiv \int_0^s \tau_r(x)dr.
\end{equation}
Considering the backward/forward differential equation 
\begin{equation}
\label{EvolJump}
\dot{f}_s = -\emph{{\bf J}}_s f_s, \;0\leq s\leq t\leq T, \quad f_t \textrm{ given,}
\end{equation}
as a perturbation of the integrable equation $\dot{f}_s = \tau_s f_s$, one sees that equation \eqref{EvolJump} is formally equivalent to the integral equation
\begin{equation}
\label{IntegralEvolJump}
f_s = e^{T_s-T_t}f_t + \int_s^t e^{T_s-T_r}\emph{{\bf L}}_rf_r\,dr.
\end{equation}

Given some positive function $h$, set $B_h = \bigl\{f\,;\,\sup\frac{|f|}{h}<\infty\bigr\}$, and define $\|f\|_h = \sup\frac{|f|}{h}$, for $f\in B_h$. The space $\bigl(B_h,\|.\|_h\bigr)$ is a Banach space. Define also $B_h^0$ as the set of functions $f\in B_h$ such that $\frac{f}{h}$ goes to $0$ as $h$ goes to infinity. The following two theorems are part of the folkore; they are stated under this form in the appendix of Kolokoltsov's article \cite{Kolokoltsov2}.

\begin{thm}[Existence of propagators, first part]
\label{ThmPropagatorsKol1}
\textbf{1)} Suppose that there exists two continuous positive functions $h$ and $h'$, and positive constants $c$ and $c'$ such that 
  \begin{itemize}
     \item[\textit{\textbf{a.}}] $0<h'\leq h$, $h'\in B_h$,
     \item[ ]                    $\forall\,s\in [0,T], \quad \emph{{\bf J}}_s h' \leq c' h', \;\,\emph{{\bf J}}_s h \leq c h$.
  \end{itemize}
Then, given $t \in [0,T]$ and some function $u_t \in B_h$, the minimal solution of the \textit{backwards/forwards} integral problem \eqref{IntegralEvolJump} with final/initial condition $u_t$ is of the form $\{S_{s,t}u_t\}_s$ for some bounded operators $S_{s,t}$ on $(B_h^0,\|.\|_h)$ depending \textit{continuously} on $s$ and $t$, with norm no greater than $e^{c|t-s|}$.
\end{thm}

If now one considers the backward/forward differential equation 
$$
\dot{f}_s = -\Lambda_s f_s = -\bigl(\emph{{\bf J}}_s-\emph{\bf M}_s\bigr)f_s, \;0\leq s\leq t\leq T, \quad f_t=f \;\textrm{ given,}
$$
as a perturbation of equation \eqref{EvolJump}, the preceding differential equation is formally equivalent to the integral equation
\begin{equation}
\label{IntEqEvol}
f_s = S_{s,t}f - \int_s^t S_{s,r}\emph{\bf M}_rf_r\,dr.
\end{equation}

\begin{thm}[Existence of propagators, second part]
\label{ThmPropagatorsKol2}
\textbf{\textit{2)}} Suppose, in addition to the hypothesis of theorem \ref{ThmPropagatorsKol1}, that the following hypothesis on the perturbations $\emph{\bf M}_s$ hold. 
  \begin{itemize}
     \item[\textit{\textbf{b.}}] The family $\{\emph{\bf M}_s\}_{0 \leq s \leq T}$ is a bounded family of linear transforms of $(B_h,\|.\|_h)$.
  \end{itemize}
Denote by $\|\emph{\bf M}_s\|_h$ the norm operator of $\emph{\bf M}_s$. Then the series
$$
U_{s,t}f = S_{s,t}f -\int_s^t S_{s,r}\emph{\bf M}_r S_{r,t}f\,dr + \int_{s\leq r_1 \leq r_2 \leq t} S_{s,r_1}\emph{\bf M}_{r_1}S_{r_1,r_2}\emph{\bf M}_{r_2}S_{r_2,t}f\,dr_1dr_2 + \cdots,
$$
converges in $(B_h,\|.\|_h)$ for any $f\in B_h$. It defines a propagator on $\left(B_h^0,\|.\|_{h}\right)$ depending \textit{continuously} on $s$ and $t$, and with norm 
\begin{equation}
\label{UpperBound}
\leq e^{\bigl(c + \underset{s \leq r \leq t}{\sup} \|\emph{\bf M}_r\|_h\bigr)|t-s|}
\end{equation}
The map $s\mapsto U_{s,t}f$ is the minimal solution of the \textit{backwards/forwards} integral problem \eqref{IntEqEvol} with final/initial condition $f$.

\textbf{\textit{3)}} If finally
  \begin{itemize}
     \item[\textbf{\textit{c}}.] $\bullet$ for any $f\in B_{h'}$, the function $s \mapsto \emph{\bf J}_sf \in (B_h^0,\|.\|_h)$ is well defined and continuous,
     \item[ ] $\bullet$ each $\emph{\bf M}_s$ sends continuously $(B_h,\|.\|_h)$ in $(B_{h'},\|.\|_{h'})$,
  \end{itemize}
\noindent then for any $f \in B_{h'}$, the function $s \mapsto U_{s,t}f \in \left(B_h^0,\|.\|_{h}\right)$ is differentiable, with derivative $-\Lambda_sU^{s,t}f$. It is also differentiable as a function of $t$, with derivative $U_{s,t}\Lambda_tf$.
\end{thm}

\bigskip

\noindent \textbf{Acknoledgements.} The author would like to thank James Norris for his interest about this work.

\end{document}